\documentclass[11pt]{amsart}

\usepackage{amsmath} 
\usepackage{amssymb}

\newtheorem{theorem}{Theorem}[section] 
\newtheorem{claim}{Claim}[theorem]
\newtheorem{mclaim}{Main Claim}[theorem]
 
\newtheorem{proposition}[theorem]{Proposition}

\theoremstyle{definition}
\newtheorem{definition}[theorem]{Definition}

\newtheorem{problem}{Problem}[section]
\newtheorem{question}[theorem]{Question}

\theoremstyle{remark}
\newtheorem{remark}[theorem]{Remark}
\newtheorem{discussion}[theorem]{Discussion}
\newtheorem{conclusion}[theorem]{Conclusion}

\newtheorem{hypothesis}[theorem]{Hypothesis}

\numberwithin{equation}{section}

\newcommand{\hd}{{\rm hd}}
\newcommand{\hL}{{\rm hL}}
\newcommand{\Aut}{{\rm Aut}}
\newcommand{\End}{{\rm End}}
\newcommand{\Inc}{{\rm Inc}}
\newcommand{\Ded}{{\rm Ded}}
\newcommand{\bd}{{\rm bd}}
\newcommand{\Ens}{{\rm Ens}}
\newcommand{\cf}{{\rm cf}}
\newcommand{\bB}{{\mathbb B}}
\newcommand{\bU}{{\bf U}}
\newcommand{\bJ}{{\bf J}}
\newcommand{\pcf}{{\rm pcf}}
\newcommand{\Atom}{{\rm Atom}}
\newcommand{\pp}{{\rm pp}}
\newcommand{\Reg}{{\rm Reg}}
\newcommand{\cov}{{\rm cov}}
\newcommand{\Dom}{{\rm Dom}}
\newcommand{\Rang}{{\rm Rang}}
\newcommand{\level}{{\rm level}}
\newcommand{\bt}{{\bf t}}
\newcommand{\conc}{{}^\frown\!}
\newcommand{\tcf}{{\rm tcf}}
\newcommand{\llg}{{\ell g}}
\newcommand{\acc}{{\rm acc}}

\newcommand{\moj}{\genfrac{}{}{0pt}1}

\title{Constructing Boolean Algebras for cardinal invariants}
\author{Saharon Shelah}
\address{Institute of Mathematics\\
 The Hebrew University of Jerusalem\\
 Jerusalem 91904, Israel\\
 and  Department of Mathematics\\
 Rutgers University\\
 New Brunswick, NJ 08854, USA}
\email{shelah@math.huji.ac.il}
\urladdr{http://www.math.rutgers.edu/$\sim$shelah}
\thanks{I would like to thank Alice Leonhardt for the beautiful typing.\\ 
This research was partially supported by the Israel Science Foundation
founded by the Israel Academy of Sciences and Humanities. Publication 641}  


\begin{document}
\begin{abstract}
We construct Boolean Algebras answering some questions of J.~Donald Monk on
cardinal invariants. The results are proved in ZFC (rather than giving
consistency results).  We deal with the existence of superatomic Boolean 
Algebras with ``few automorphisms'', with entangled sequences of linear orders,
and with semi-ZFC examples of the non-attainment of the spread (and
$\hL,\hd$).
\end{abstract}
\maketitle

\section*{Annotated Content}
\begin{enumerate}
\item[\S1] {\em A superatomic Boolean Algebra with fewer automorphisms
than endomorphisms.}

\noindent We prove in ZFC that for some superatomic Boolean Algebra $\bB$
we have $\Aut(\bB)<\End(\bB)$.  This solves \cite[Problem 76, p.291]{M2} of
Monk. 

\item[\S2] {\em A superatomic Boolean Algebra with fewer automorphisms than
elements}

\noindent We prove in ZFC that for some superatomic Boolean Algebra $\bB$, 
we have $\Aut(\bB)<|\bB|$. This solves \cite[Problem 80, p.291]{M2} of Monk.

\item[\S3] {\em On entangledness}

\noindent We prove that if $\mu<\kappa\le\chi<\Ded(\mu)$ and $2^\mu<
\lambda$, and $\kappa$ is regular, and $\lambda\leq\bU_{J^{\bd}_\kappa}(
\chi)$ (see Definition \ref{3.2}), then $\Ens(\kappa,\lambda)$, i.e., there
is an entangled sequence of $\lambda$ linear orders each of cardinality
$\kappa$. The reader may think of the case 
\[\mu=\aleph_0,\quad\kappa=\cf(\chi)<\chi=2^\mu=2^{<\kappa}<2^\kappa,\quad
\mbox{ and }\quad\lambda=\chi^+.\]
Note that the existence of entangled linear orders is connected to the
problem whether always $\prod\limits_{i<\theta}\Inc(\bB_i)/D \geq\Inc(
\prod\limits_{i<\theta}\bB_i/D)$ for an ultrafilter $D$ on $\theta$. We rely
on quotations of some pcf results.

\item[\S4] {\em On attainment of spread}

\noindent We construct Boolean Algebras with the spread not obtained under
ZFC + ``GCH is violated strongly enough, even just for regular cardinals'';
so the consistency strength is ZFC. We consider this a semi--ZFC answer.
\end{enumerate}

\section{A superatomic Boolean Algebra with fewer automorphisms than
endomorphisms} 

Rubin has proved that if $\diamondsuit_{\lambda^+}$, then there is a
superatomic Boolean algebra with few automorphisms. We give here a
construction in ZFC.  

We use some notions of \cite{Sh:g}, they can be found in \cite{Sh:E12}; in
particular $\bJ_\theta[{\mathfrak a}]=\bJ_{<\theta}[{\mathfrak a}]+
({\mathfrak a}\setminus {\mathfrak b}_\theta[{\mathfrak a}])$. For this
section we assume  

\begin{hypothesis}
\label{1.1} 
\begin{enumerate}
\item[(a)] $\bar{\lambda}=\langle\lambda_i:i<\delta\rangle$ is a strictly
increasing sequence of regular cardinals larger than $\delta$; let
${\mathfrak a}=\{\lambda_i:i<\delta\}$.
\item[(b)] $\lambda_0>2^{|\delta|}$, or at least $\lambda_0>|\pcf({\mathfrak
a})|$.  
\end{enumerate}
\end{hypothesis}

The main combinatorial point of our construction is given by the following
observation. 

\begin{proposition}
\label{1.2}
There are sequences $\langle\bar{f}^\theta:\theta\in\pcf({\mathfrak a})
\rangle$ and $\langle{\mathfrak b}_\theta[{\mathfrak a}]:\theta\in
\pcf({\mathfrak a})\rangle$ such that
\begin{enumerate}
\item[(a)] $\bar{f}^\theta=\langle f^\theta_\alpha:\alpha<\theta\rangle
\subseteq\prod{\mathfrak a}$ is a $<_{\bJ_{\theta}[{\mathfrak
a}]}$--increasing cofinal sequence, $\langle{\mathfrak b}_\theta[{\mathfrak
a}]:\theta\in\pcf({\mathfrak a})\rangle$ is a generating sequence,
\item[(b)] $f^\theta_\alpha\restriction ({\mathfrak a}\setminus {\mathfrak
b}_\theta[{\mathfrak a}])$ is constantly zero, 
\item[(c)] if $\theta_1<\theta_2$, $\alpha_2<\theta_2$ then
\[f^{\theta_2}_{\alpha_2}\restriction ({\mathfrak b}_{\theta_1}[{\mathfrak
a}]\cap {\mathfrak b}_{\theta_2}[{\mathfrak a}])\in\{f^{\theta_1}_{\alpha_1} 
\restriction ({\mathfrak b}_{\theta_1}[{\mathfrak a}]\cap {\mathfrak b}_{
\theta_2}[{\mathfrak a}]):\alpha_1<\theta_1\}.\]
\item[(d)] for $\theta\in\pcf({\mathfrak a})$ and $\lambda\in{\mathfrak
b}_\theta[{\mathfrak a}]$, $f^\theta_\alpha(\lambda)$ is a limit ordinal $>
\sup(\lambda\cap {\mathfrak a})$. 
\item[(e)] if $\theta_1<\theta_2$, both in $\pcf({\mathfrak a})$, then there
are $n<\omega$, $\sigma_1,\ldots,\sigma_n\leq\theta$ (all from
$\pcf({\mathfrak a})$) such that ${\mathfrak b}_{\theta_1}[{\mathfrak a}]
\cap {\mathfrak b}_{\theta_2}[{\mathfrak a}]=\bigcup\limits_{k=1}^n
{\mathfrak b}_{\sigma_k}[{\mathfrak a}]$.  
\end{enumerate}
\end{proposition}

\begin{proof}  
Let ${\mathfrak a}'=\pcf({\mathfrak a})$, so $|{\mathfrak a}'|<\min({
\mathfrak a}')$ and $\pcf({\mathfrak a}')={\mathfrak a}'$ (by \cite[Ch.I,
1.11]{Sh:g}). We can find a generating sequence $\langle{\mathfrak b}_\theta
[{\mathfrak a}']:\theta\in {\mathfrak a}'\rangle$ (by \cite[Ch.VIII,
2.6]{Sh:g}), and hence a closed smooth one (by \cite[Ch.I,
3.8(3)]{Sh:g}). Now repeat the proof of \cite[Ch.II, 3.5]{Sh:g} or see
\cite{Sh:E12}.  Note that ``smooth'' means 
\[\sigma\in {\mathfrak b}_\theta[{\mathfrak a}']\quad \Rightarrow\quad
{\mathfrak b}_\sigma[{\mathfrak a}']\subseteq {\mathfrak
b}_\theta[{\mathfrak a}'],\]
``closed'' means ${\mathfrak b}_\theta[{\mathfrak a}]=\pcf({\mathfrak
b}_\theta[{\mathfrak a}])\cap {\mathfrak a}$; together clause (e) follows.
\end{proof}

\begin{definition}
\label{1.3}
Let $\langle\bar{f}^\theta:\theta\in\pcf({\mathfrak a})\rangle$ and
$\langle{\mathfrak b}_\theta[{\mathfrak a}]:\theta\in\pcf({\mathfrak a})
\rangle$ be sequences given by \ref{1.2} (so they satisfy the demands
(a)--(e) there). 
\begin{enumerate}
\item For $\ell\in \{0,1\}$, $\theta\in\pcf({\mathfrak a})$ and $\alpha\leq
\theta$ we define the Boolean ring ${\mathcal B}^\ell_{\theta,\alpha}$ of
subsets of $\sup({\mathfrak a})$. We do this by induction on $\theta$, and
for each $\theta$ by induction on $\alpha$ as follows. 
\begin{enumerate}
\item[(a)] If $\theta=\min({\mathfrak a})$, $\alpha=0$, then ${\mathcal
B}^\ell_{\theta,\alpha}$ is the Boolean ring (of subsets of $\sup({\mathfrak
a})$) generated by 
\[\{[\sup({\mathfrak a}\cap\lambda),\gamma):\lambda\in {\mathfrak a},\
\sup({\mathfrak a}\cap\lambda)<\gamma<\lambda\},\] 
that is the closure of the above family under $x\cap y$, $x\cup y$, $x-y$.
\item[(b)] If $\theta\in\pcf({\mathfrak a})\setminus\{\min({\mathfrak
a})\}$, then let ${\mathcal B}^\ell_{\theta,0}=\bigcup\limits_{\sigma\in
\theta \cap\pcf({\mathfrak a})}{\mathcal B}^\ell_{\sigma,\sigma}$. 
\item[(c)] If $\theta\in\pcf({\mathfrak a})$, $\alpha<\theta$ is a limit
ordinal then ${\mathcal B}^\ell_{\theta,\alpha}=\bigcup\limits_{\beta<
\alpha}{\mathcal B}^\ell_{\theta,\beta}$. 
\item[(d)] If $\theta\in\pcf({\mathfrak a})$, $\alpha=\theta$, then
${\mathcal B}^\ell_{\theta,\alpha}$ is the Boolean ring generated by 
\[\bigcup\limits_{\beta<\alpha}{\mathcal B}^\ell_{\theta,\beta}\cup\big\{
\bigcup\limits_{\lambda\in {\mathfrak b}_\theta[{\mathfrak a}]}[\sup(\lambda
\cap {\mathfrak a}),\lambda)\big\}.\]
\item[(e)] If $\theta\in\pcf({\mathfrak a})$, $\alpha=\beta+1<\theta$, then
\begin{enumerate}
\item[(i)]  ${\mathcal B}^0_{\theta,\alpha}$ is the Boolean ring generated by 
\[{\mathcal B}^1_{\theta,\beta}\cup\big\{\bigcup\limits_{\lambda\in
{\mathfrak b}_\theta[{\mathfrak a}]}[\sup({\mathfrak a}\cap\lambda),
f^\theta_\beta(\lambda))\big\},\]
\item[(ii)] ${\mathcal B}^1_{\theta,\alpha}$ is the Boolean ring generated
by 
\[{\mathcal B}^0_{\theta,\alpha}\cup\big\{\bigcup\limits_{\lambda\in
{\mathfrak b}_\theta[{\mathfrak a}]}[\sup({\mathfrak a}\cap\lambda),
f^\theta_\beta(\lambda)+1)\big\}.\] 
\end{enumerate}
\end{enumerate}
\item We let
\[\begin{array}{lcll}
x^0_{\theta,\beta}&=&\bigcup\limits_{\lambda\in {\mathfrak b}_\theta[{
\mathfrak a}]}[\sup({\mathfrak a}\cap\lambda),f^\theta_\beta(\lambda))&
\mbox{ for }\theta\in\pcf({\mathfrak a}),\ \beta<\theta,\\ 
x^1_{\theta,\beta}&=&\bigcup\limits_{\lambda\in {\mathfrak b}_\theta[{
\mathfrak a}]}[\sup({\mathfrak a}\cap\lambda),f^\theta_\beta(\lambda)+1)&
\mbox{ for }\theta\in\pcf({\mathfrak a}),\ \beta<\theta,\\ 
y_\theta&=&\bigcup\limits_{\lambda\in {\mathfrak b}_\theta[{\mathfrak a}]}
[\sup(\lambda\cap{\mathfrak a}),\lambda)&\mbox{ for }\theta\in\pcf({
\mathfrak a}),\\
z_\alpha&=&[\sup(\lambda\cap{\mathfrak a}),\alpha)&\mbox{ for }\alpha\in
  [\sup(\lambda\cap {\mathfrak a}),\lambda),\ \lambda\in {\mathfrak a}.
  \end{array}\]
\item $\bB^\ell_{\theta,\alpha}$ is the Boolean algebra of subsets of
$\sup({\mathfrak a})$ generated by ${\mathcal B}^\ell_{\theta,\alpha}$, and
$\bB^\ell$ stands for the Boolean Algebra of subsets of $\sup({\mathfrak
a})$ generated by ${\mathcal B}^\ell_{\max\pcf({\mathfrak a}),
\max\pcf({\mathfrak a})}$.  

[After we shall note that $\bB^0=\bB^1$ (in \ref{1.5}) we can write $\bB^0=
\bB=\bB^1$.] 
\end{enumerate}
\end{definition}

\begin{proposition}
\label{1.5} 
\begin{enumerate}
\item ${\mathcal B}^\ell_{\theta,\theta}$ is increasing in $\theta$, and for
a fixed $\theta$, ${\mathcal B}^\ell_{\theta,\alpha}$ is increasing in
$\alpha$ and is actually a Boolean ring of subsets of $\sup({\mathfrak
a})$. 
\item ${\mathcal B}^m_{\theta,\alpha}$ is the Boolean ring generated by
\[\begin{array}{l}
\{y_\sigma:\sigma\in\pcf({\mathfrak a})\cap\theta\ \mbox{ or }\ \sigma=
\alpha=\theta\}\cup\{z_\alpha:\alpha<\sup({\mathfrak a})\}\cup{}\\ 
\{x^\ell_{\sigma,\alpha}:\sigma<\theta,\ \sigma\in\pcf({\mathfrak a}),\
\alpha<\sigma,\ \ell<2\}\cup{}\\ 
\{x^\ell_{\theta,\beta}:\beta+1<\alpha\ \&\ \ell<2\mbox{ or }\beta+1=\alpha\
\&\ \ell\leq m\}. 
  \end{array}\]
\item If $\alpha$ is zero or limit $\leq\theta\in\pcf({\mathfrak a})$, then
${\mathcal B}^0_{\theta,\alpha}={\mathcal B}^1_{\theta,\alpha}$ and $\bB^0_{
\theta,\alpha}=\bB^1_{\theta,\alpha}$.
\item If $(\theta_1,\alpha_1,\ell_1)\leq_{\ell ex}(\theta_2,\alpha_2,
\ell_2)$ and $a_i\in {\mathcal B}^{\ell_i}_{\theta_i,\alpha_i}$ for $i
=1,2$, {\em then\/} $a_1\cap a_2\in {\mathcal B}^{\ell_1}_{\theta_1,
\alpha_1}$.
\item If $\ell_i\in\{0,1\}$, $\alpha_i\leq\theta_i\in\pcf({\mathfrak a})$
for $i=1,2$ and  
\[\theta_1<\theta_2\vee(\theta_1=\theta_2\ \&\ \alpha_1<\alpha_2)\vee
(\theta_1=\theta_2\ \&\ \alpha_1=\alpha_2\ \&\ \ell_1<\ell_2),\]
{\em then\/} ${\mathcal B}^{\ell_1}_{\theta_1,\alpha_1}$ is an ideal of
${\mathcal B}^{\ell_2}_{\theta_2,\alpha_2}$. 
\end{enumerate}
\end{proposition}

\begin{proof} (1)--(3)\qquad Straightforward.
\medskip

\noindent (4)\qquad First note that it is enough to show the assertion under
an additional demand that $a_1,a_2$ are among the generators of the Boolean
rings ${\mathcal B}^{\ell_1}_{\theta_1,\alpha_1},{\mathcal B}^{\ell_2}_{
\theta_2,\alpha_2}$, respectively, as listed in part (2).
\smallskip

\noindent{\sc Case 1:}\qquad One of $a_1,a_2$ is $z_\alpha$ for some
$\alpha<\sup({\mathfrak a})$.\\
Then the other is either $y_\theta$, or $z_\beta$, or $x^m_{\theta,\beta}$,
and in all cases the intersection $a_1\cap a_2$ is either empty or it is
$z_{\alpha'}$ for some $\alpha'\leq\alpha$. Hence $a_1\cap a_2\in {\mathcal
B}^{\ell_1}_{\theta_1,\alpha_1}$. 
\smallskip

\noindent{\sc Case 2:}\qquad $a_1=y_{\theta'}$, $a_2=y_{\theta''}$ for some
$\theta',\theta''\in\pcf({\mathfrak a})$.\\
If $\theta''\leq\theta'$ then, as $a_1\in {\mathcal B}^{\ell_1}_{\theta_1,
\alpha_1}$, we easily get $a_2\in {\mathcal B}^{\ell_1}_{\theta_1,\alpha_1}$
and thus the intersection $a_1\cap a_2$ is in this Boolean ring.

So we may assume that $\theta'<\theta''$. It follows from \ref{1.2}(e) that
there are $\sigma_1,\ldots,\sigma_n\leq\theta'$ such that 
\[{\mathfrak b}_{\theta'}[{\mathfrak a}]\cap {\mathfrak b}_{\theta''}[
{\mathfrak a}]=\bigcup_{k=1}^n {\mathfrak b}_{\sigma_k}[{\mathfrak a}].\]
Then $a_1\cap a_2=y_{\sigma_1}\cup\ldots\cup y_{\sigma_n}$ and
$y_{\sigma_1},\ldots,y_{\sigma_n}\in {\mathcal B}^\ell_{\theta_1,\alpha_1}$, 
so we are done.
\smallskip

\noindent{\sc Case 3:}\qquad $a_1=y_{\theta'}$, $a_2=x^m_{\theta'',\beta}$
for some $\theta',\theta''\in\pcf({\mathfrak a})$, $m<2$,
$\beta<\theta''$.\\
If $\theta''\leq\theta'$ then $a_2\in{\mathcal B}^{\ell_1}_{\theta_1,
\alpha_1}$ and we are done; so assume that $\theta'<\theta''$. It follows
from \ref{1.2}(c) that then $f^{\theta''}_\beta\restriction (b_{\theta'}[
{\mathfrak a}]\cap b_{\theta''}[{\mathfrak a}])=f^{\theta'}_\alpha
\restriction (b_{\theta'}[{\mathfrak a}]\cap b_{\theta''}[{\mathfrak a}])$
for some $\alpha<\theta'$. Like in Case 2, one shows that  $y_{\theta'}\cap
y_{\theta''}\in {\mathcal B}^{\ell_1}_{\theta_1,\alpha_1}$; also
$x^m_{\theta',\alpha}\in {\mathcal B}^{\ell_1}_{\theta_1,\alpha_1}$. But 
now $a_1\cap a_2=x^m_{\theta',\alpha}\cap (y_{\theta'}\cap y_{\theta''})\in
{\mathcal B}^{\ell_1}_{\theta_1,\alpha_1}$. 
\smallskip

\noindent{\sc Case 4:}\qquad $a_1=x^m_{\theta',\beta}$, $a_2=y_{\theta''}$
for some $\theta',\theta''\in\pcf({\mathfrak a})$, $\beta<\theta'$, $m<2$.\\
If $\theta''<\theta'$ then $a_2\in {\mathcal B}^{\ell_1}_{\theta_1,
\alpha_1}$, and if $\theta''=\theta'$ then $a_1\cap a_2=a_1$. So we may
assume $\theta'<\theta''$. If ${\mathfrak b}_{\theta'}[{\mathfrak a}]
\subseteq {\mathfrak b}_{\theta''}[{\mathfrak a}]$, then clearly $a_1\cap
a_2=a_1$ and we are done, so suppose otherwise. Then, using \ref{1.2}(e), we
find $\sigma_1,\ldots,\sigma_n<\theta'$ such that ${\mathfrak b}_{\theta'}[
{\mathfrak a}]\cap {\mathfrak b}_{\theta''}[{\mathfrak a}]=
\bigcup\limits_{k=1}^n {\mathfrak b}_{\sigma_k}[{\mathfrak a}]$. Since, in
this case, all $\sigma_k$ are smaller than $\theta'$ and $x^m_{\theta',
\beta}\cap y_{\theta''}=\bigcup\limits_{k=1}^n y_{\sigma_k}\cap
x^m_{\theta',\beta}$, we easily conclude $a_1\cap a_2\in {\mathcal
B}^{\ell_1}_{\theta_1,\alpha_1}$. 
\smallskip

\noindent{\sc Case 5:}\qquad $a_1=x^{m'}_{\theta',\beta'}$, $a_2=x^{m''}_{
\theta'',\beta''}$ for some $\theta',\theta''\in\pcf({\mathfrak a})$,
$\beta'<\theta'$, $\beta''<\theta''$ and $m',m''<2$.\\
If $(\theta'',\beta'',m'')\leq_{\ell ex} (\theta',\beta',m')$ then we are
easily done.

If $\theta''=\theta'$, $\beta''=\beta'$ and  $0=m'<m''=1$, then clearly
$a_1\cap a_2=a_1$.

Assume that $\theta''=\theta'$, $\beta'<\beta''$. Then, by \ref{1.2}(a), we 
find $\mu_0,\ldots,\mu_{k-1}\in\theta'\cap\pcf({\mathfrak a})$ such that 
\[\{\mu\in{\mathfrak a}:f^{\theta''}_{\beta''}(\mu)\leq f^{\theta'}_{
\beta'}(\mu)\}\subseteq\bigcup_{j<k}{\mathfrak b}_{\mu_j}[{\mathfrak a}]\cup
({\mathfrak a}\setminus{\mathfrak b}_{\theta'}).\] 
Then clearly
\[a_1\cap a_2= \big(\bigcup\limits_{j<k}a_1\cap a_2\cap y_{\mu_j}
\big)\cup \big(a_1\setminus\bigcup\limits_{j<k} y_{\mu_j}\big).\] 
Also, for $j<k$, we have 
\[y_{\mu_j}\in {\mathcal B}^{\ell_1}_{\theta_1,\alpha_1}\quad\mbox{ and
}\quad a_1\cap a_2\cap y_{\mu_j}=(a_1\cap y_{\mu_j})\cap (a_2\cap
y_{\mu_j}),\] 
and the sets $a_1\cap y_{\mu_j}$ and $a_2\cap y_{\mu_j}$ are in ${\mathcal
B}^{\ell_1}_{\theta_1,\alpha_1}$ by (suitably applied) case 3. So we can
easily finish.

The only remaining possibility is $\theta'<\theta''$. By \ref{1.2}(a) we may
pick $\gamma<\theta'$ such that 
\[f^{\theta''}_{\beta''}\restriction ({\mathfrak b}_{\theta''}[{\mathfrak a}]
\cap {\mathfrak b}_{\theta'}[{\mathfrak a}])=f^{\theta'}_\gamma\restriction
({\mathfrak b}_{\theta''}[{\mathfrak a}]\cap {\mathfrak b}_{\theta'}[
{\mathfrak a}]).\]
Then $a_1\cap a_2=x^{m'}_{\theta',\beta'}\cap x^{m''}_{\theta',\gamma}\cap
y_{\theta'}\cap y_{\theta''}$. By the discussion above we know that
$x^{m'}_{\theta',\beta'}\cap x^{m''}_{\theta',\gamma}\in {\mathcal
B}^{\ell_1}_{\theta_1,\alpha_1}$. Now, if $y_{\theta'}\subseteq
y_{\theta''}$ then $a_1\cap a_2= x^{m'}_{\theta',\beta'}\cap
x^{m''}_{\theta',\gamma}$ and we are done. Otherwise, $y_{\theta'}\cap
y_{\theta''}\in {\mathcal B}^0_{\theta',0}\subseteq {\mathcal B}^{\ell_1}_{
\theta_1,\alpha_1}$ (compare Case 2), and again we easily get the required
conclusion.  
\medskip

\noindent (5)\qquad Follows. 
\end{proof}

\begin{proposition}
\label{1.6} \begin{enumerate}
\item ${\mathcal B}^\ell_{\theta,\alpha}$ is a superatomic Boolean ring with
$\{\{\gamma\}:\gamma<\sup({\mathfrak a})\}$ as the set of atoms.
\item $\bB^\ell_{\theta,\alpha}$ is a superatomic Boolean algebra, in
particular $\bB^\ell$ is.
\item If $\alpha,\beta<\theta\in\pcf({\mathfrak a})$ and $\gamma=
\omega^\beta$ (ordinal exponentiation; so $\gamma<\theta$ and
$\alpha+\gamma<\theta$), {\em then\/} the rank of $x^0_{\theta,\alpha+
\gamma}-x^0_{\theta,\alpha}$ is $\geq\beta$. 
\end{enumerate}
\end{proposition}

\begin{proof}  1)\quad Straight by induction on $\theta$ and for a fixed
$\theta$ by induction on $\alpha\leq\theta$ using \ref{1.5}(5).
\medskip

\noindent 2)\quad Follows. 
\medskip

\noindent 3)\quad Easy by induction on $\beta$. 
\end{proof}

\begin{proposition}
\label{1.7}
\begin{enumerate}
\item The algebra $\bB$ has exactly $\sup({\mathfrak a})$ atoms, so 
\[|\Atom(\bB^\ell)|=\sup({\mathfrak a}).\]
\item $|\bB|=\max\pcf({\mathfrak a})$.
\item $|\Aut(\bB)|\le 2^{\sup({\mathfrak a})}$.
\end{enumerate}
\end{proposition}

\begin{proof} Parts 1), 2) should be clear. Part 3) holds as the algebra
$\bB$ has $\sup({\mathfrak a})$ atoms by part (1) (and two distinct
automorphisms of $\bB$ differ on an atom). 
\end{proof}

\begin{proposition}
\label{1.8}  
The algebra $\bB$ has $2^{\max\pcf({\mathfrak a})}$ endomorphisms.
\end{proposition}

\begin{proof}
Let $Z\subseteq\max\pcf({\mathfrak a})$. We define an endomorphism $T_Z\in
\End(\bB)$ by describing how it acts on the generators. We let:
\[\begin{array}{lcl}
T_Z(z_\alpha)&=& z_\beta\quad\mbox{ if $\beta$ is maximal such that }\beta
\leq\alpha<\beta+\omega\mbox{ and:}\\
&&\beta=0\mbox{ or }\beta\mbox{ limit or }\alpha=\beta\in
\bigcup\limits_{\lambda\in {\mathfrak a}}[\sup({\mathfrak a}\cap\lambda),
\sup({\mathfrak a}\cap\lambda)\!+\!\omega),\\
T_Z(y_\theta)&=&y_\theta,\\
T_Z(x^0_{\theta,\alpha})&=&x^0_{\theta,\alpha},\\
&&\\
T_Z(x^1_{\theta,\alpha})&=&\left\{
                  \begin{array}{ll}
            x^0_{\theta,\alpha}&\mbox{if }\theta<\max\pcf({\mathfrak a}),\\ 
            x^0_{\theta,\alpha}&\mbox{if }\theta=\max\pcf({\mathfrak a}),\
                                 \alpha\notin Z,\\ 
            x^1_{\theta,\alpha}&\mbox{if }\theta=\max\pcf({\mathfrak a}),\
\alpha\in Z. 
		  \end{array}\right.
  \end{array}\]
One easily checks that the above formulas correctly define an element of
$\End(\bB)$. Clearly $Z_1\neq Z_2$ implies $T_{Z_1}\neq T_{Z_2}$ and we are
done. 
\end{proof}

So we can answer (in ZFC) Monk's question \cite[Problem 76, pages 259,
291]{M2}. 

\begin{conclusion}
\label{1.9}
Assume that $\mu$ is a strong limit singular cardinal, and $\cf(\mu)>
\aleph_0$ (or just $\pp^+(\mu)=(2^\mu)^+$, so most of those with $\cf(\mu)=
\aleph_0$ are OK) and $\mu<\kappa=\cf(\kappa)\leq 2^\mu<2^\kappa$ (always
such $\mu$ exists and for each such $\mu$ such $\kappa$ exists). {\em
Then\/} there is a superatomic Boolean Algebra $\bB$ such that: 
\begin{enumerate}
\item[(a)] $|\bB|=\kappa$,
\item[(b)] $|\Atom(\bB)|=\mu$,
\item[(c)] $|\Aut(\bB)|\leq 2^\mu$,
\item[(d)] $|\End(\bB)|=2^\kappa$.
\end{enumerate}
\end{conclusion}

\begin{proof} 
We can find ${\mathfrak a}\subseteq\Reg\cap\mu$ such that $|{\mathfrak a}|=
\cf(\mu)$ and $\kappa=\max\pcf({\mathfrak a})$. Why? We know 
\[2^\mu=\mu^{\cf(\mu)}=\cov(\mu,\mu,(\cf(\mu))^+,2)\ge\cov(\mu,\mu,\cf(\mu)^+, 
\cf(\mu)),\]
and now we use \cite[Ch.II, 5.4]{Sh:g} when $\cf(\mu)>\aleph_0$;
see \cite[6.5]{Sh:E12} for references on the $\cf(\mu)=\aleph_0$ case). 
\end{proof}

\section{A superatomic Boolean Algebra with fewer automorphisms than
elements}

Monk has asked (\cite[Problem 80, p.291,260]{M2}) if there may be a
superatomic Boolean Algebra $|\bB|$ with ``few'' (i.e., $<|\bB|$)
automorphisms. Remember that $\Aut(\bB)\geq|\Atom(\bB)|$ if $|\Atom(\bB)|\ne
1$. 

In this section we answer this question by showing that, in ZFC, there is a
superatomic Boolean Algebra $\bB$ with $\Aut(\bB)<|\bB|$. Moreover, there
are such Boolean Algebras in many cardinals. 

For our construction we assume the following:
\begin{hypothesis}
\label{2.1A} 
\begin{enumerate}
\item[$(\alpha)$] $\mu$ is a strong limit singular cardinal of cofinality
$\aleph_0$, 
\item[$(\beta)$]  $\lambda=2^\mu$, $\kappa\leq\lambda$,
\item[$(\gamma)$] $T$ is a tree with $\kappa$ levels, $\leq\lambda$ nodes
and the number of its $\kappa$-branches is $\chi>\lambda$, and $T$ has a
root. 
\end{enumerate}
\end{hypothesis}

Note that there are many $\mu$ as in clause $(\alpha)$ of \ref{2.1A}, and
then we can choose $\lambda=2^\mu$ and, e.g., $\kappa=\min\{\kappa:2^\kappa
>\lambda\}$, $T= {}^{\kappa>}2$.

\begin{theorem}
\label{2.1B}
There is a superatomic Boolean Algebra $\bB$ such that:
\[|\bB|=\chi\quad\mbox{ and }\quad |\Atom(\bB)|=|T|+\mu\le |\Aut(\bB)|\le
\lambda.\]
\end{theorem}

\begin{proof}
Let $T^+=T\cup\lim_\kappa(T)$, so $|T^+|=\chi$. Let
\[\begin{array}{lr}
{\mathcal F}=\bigl\{f:&f\mbox{ is a one-to-one function, }\Dom(f)\subseteq
T \times\mu,\quad\\ 
&|\Dom(f)|=\mu,\ \Rang(f)\subseteq T\times\mu\setminus\Dom(f)\bigr\}.
  \end{array}\]
Clearly $|{\mathcal F}|\le |T\times\mu|^\mu\leq \lambda^\mu=(2^\mu)^\mu=
2^\mu =\lambda$.

Let $x_{(t,\alpha)}=x_{t,\alpha}=\{(t,\alpha)\}$ (for $t\in T$ and $\alpha<
\mu$) and let $z_t=\{(s,\alpha):s<_T t \mbox{ and }\alpha<\mu\}$ (for $t\in
T^+$). (Note that if $t_1,t_2$ are immediate successors of $s$, then
$z_{t_1}=z_{t_2}$; also the family $\{z_t:t\in T^+\}$ is closed under
intersections.)  

\begin{claim}
\label{cl1}
There is a family ${\mathcal A}_2\subseteq [T\times \mu]^{\aleph_0}$ such
that 
\begin{enumerate}
\item[(a)] if $y',y'' \in {\mathcal A}_2$ are distinct, then $y'\cap y''$ is  
finite, 
\item[(b)] if $t\in T^+$ and $y\in {\mathcal A}_2$, then 
\[y\subseteq z_t\ \vee\ |y\cap z_t|<\aleph_0,\]
\item[(c)] if $Y\in [T\times\mu]^\mu$, and $f$ is a one-to-one function from
$Y$ to $T\times\mu\setminus Y$ (so $f\in {\mathcal F}$), {\em then} there is 
$y\in {\mathcal A}_2$ such that $f[y]$ is almost disjoint from every member
of ${\mathcal A}_2$. 
\end{enumerate}
\end{claim}

\begin{proof}[Proof of the Claim]
List ${\mathcal F}$ as $\{g_\alpha:\alpha<\lambda\}$. By induction on
$\alpha<\lambda$ we choose $y_\alpha,y'_\alpha$ such that:
\begin{enumerate}
\item[$(\alpha)$]  $y_\alpha,y'_\alpha$ are disjoint countable subsets
of $T\times\mu$,
\item[$(\beta)$]   $y_\alpha,y'_\alpha$ are almost disjoint to any $y''\in
\{y_\beta,y'_\beta:\beta<\alpha\}$, 
\item[$(\gamma)$]  $y_\alpha\subseteq\Dom(g_\alpha)$, $y'_\alpha=g_\alpha
[y_\alpha]$,
\item[$(\delta)$]  if $t\in T^+$ then either $y_\alpha\subseteq z_t$ or
$|y_\alpha\cap z_t|<\aleph_0$. 
\end{enumerate}
So assume $y_\beta,y'_\beta$ for $\beta<\alpha$ have been defined. Pick an
increasing sequence $\langle\mu_n:n<\omega\rangle$ of regular cardinals such
that $\mu=\sum\limits_{n<\omega}\mu_n$ and $2^{\mu_n}<\mu_{n+1}$.   

Choose pairwise disjoint sets $Y_n\in [\Dom(g_\alpha)]^{\mu_n}$ (for $n<
\omega$). We may replace $Y_n$ by any $Y'_n\in [Y_n] ^{\mu_n}$, and even
$Y'_n\in [Y_{k_n}]^{\mu_n}$ with a strictly increasing sequence $\langle
k_n:n<\omega\rangle$.

Let $Y_n=\{(t^n_i,\alpha^n_i):i<\mu_n\}$ be an enumeration (with no
repetitions). Without loss of generality: 
\begin{itemize}
\item the sequence $\langle\level(t^n_i):i<\mu_n\rangle$ is constant or
strictly increasing, 
\item the sequence $\langle\alpha^n_i:i<\mu_n\rangle$ is constant or
strictly increasing, and 
\item for each $n<\omega$, for some truth value $\bt_n$ we have 
\[(\forall i<j<\mu_n)(\mbox{truth value}\;(t^n_i<_T t^n_j)\equiv\bt_n).\]
\end{itemize}
[Why?  E.g.~use $\mu_{n+1}\rightarrow(\mu_n)^2_2$].  Cleaning a little more
we may demand that 
\begin{itemize}
\item for $n\neq m$, for some truth value $\bt_{m,n}$,
\[(\forall i<\mu_n)(\forall j<\mu_m)(\mbox{truth value}\;(t^n_i<_T t^m_j)=
\bt_{m,n}).\]
\end{itemize}
[Why? E.g. use polarized partition relations.]  Using Ramsey's theorem
applied to the partition $F(m,n)=\bt_{m,n}$ (and replacing $\langle\mu_n:n<
\omega\rangle$ by an $\omega$--subsequence), without loss of generality:
\[\begin{array}{rl}
\mbox{either:}&\mbox{for some }t\in T^+,\mbox{ for every }\eta\in
\prod\limits_{n<\omega}\mu_n\mbox{ we have}\\ 
&\{(t^\ell_{\eta(\ell)},\alpha^\ell_{\eta(\ell)}):\ell<\omega\}\subseteq
z_t,\\
&\\
\mbox{or:}&\mbox{for every }t\in T^+\mbox{ and }\eta\in\prod\limits_{n<
\omega}\mu_n\mbox{ we have}\\ 
&|\{(t^\ell_{\eta(\ell)},\alpha^\ell_{\eta(\ell)}:\ell<\omega)\}\cap z_t|
\le 1. 
  \end{array}\]
Next we choose $\{(t_\eta,\beta_\eta):\eta\in\prod\limits_{\ell<n}\mu_\ell
\}\subseteq Y_n$ (no repetitions) and for each $\eta\in\prod\limits_{n<
\omega} \mu_n$ we consider 
\[y_\eta=\{(t_{\eta\restriction\ell},\beta_{\eta\restriction\ell}):\ell<
\omega\},\quad y'_\eta=g_\alpha[y_\eta]\]
as candidates for $y_\alpha,y_\alpha'$, respectively. Clause $(\alpha)$
holds as $\Rang(g_\alpha)\cap\Dom(g_\alpha)=\emptyset$, clauses
$(\gamma)$ and $(\delta)$ are also trivial. So only clause $(\beta)$ may
fail. Each $\beta<\alpha$ disqualifies at most $2^{\aleph_0}$ of the
$\eta$'s, i.e., of the pairs $(y_\eta,y'_\eta)$. So only $\leq |\alpha|
\times 2^{\aleph_0}<\lambda=|{}^\omega\mu|$ of the $\eta$'s are
disqualified, so some are OK, and we are done. This finishes the proof of
the Claim.  
\end{proof}

Let ${\mathcal A}_2$ be a family given by \ref{cl1} and let 
\[{\mathcal A}_0=\{\{x\}:x \in T \times \mu\}\qquad\mbox{ and }\qquad
{\mathcal A}_1=\{z_t:t \in T^+\}\]
Our Boolean algebra $\bB$ is the Boolean Algebra of subsets of $T\times\mu$
generated by ${\mathcal A}_0\cup {\mathcal A}_1 \cup {\mathcal A}_2$.

\begin{claim}
\label{cl2}
The algebra $\bB$ is superatomic.
\end{claim}

\begin{proof}[Proof of the Claim]
Clearly, the family $I=\{b\in\bB:\bB\restriction b$ is superatomic$\}$
is an ideal in $\bB$. Plainly $x_{(t,\alpha)}\in I$ for $(t,\alpha)\in
T \times\mu$. Now, by induction on $\alpha\leq\kappa$ we prove that if
$t\in T^+$ is of level $\alpha$, then $\bB\restriction z_t$ is
superatomic.\\
If $\alpha=0$, then $z_t=\emptyset$ and this is trivial.\\
If $\alpha=\beta+1$ and $t$ is the immediate successor of $s$, then (as
$\bB\restriction z_s$ is superatomic by the induction hypothesis) it is
enough to prove that $\bB \restriction (z_t-z_s)$ is superatomic.  Now,
$\bB\restriction (z_t-z_s)$ is the Boolean Algebra of subsets of
$\{s\}\times\mu$, generated by  
\[\{\{(s,\alpha)\}:\alpha<\mu\}\cup\{y\cap(\{s\}\times\mu):y\in {\mathcal
A}_2\},\]
and we are done by \ref{cl1}(a).\\
If $\alpha$ is a limit ordinal and $\cf(\alpha)=\aleph_0$, then $(\bB
\restriction z_t)/{\rm id}_\bB(\{z_{t\restriction\beta}:\beta<\alpha\})$ is
a Boolean Algebra generated by its atoms  
\[\{y\cap z_t:y\in {\mathcal A}_2\mbox{ and } y\leq z_t\ \&\ \bigwedge_{s<t}
\neg y\leq z_s\}\]  
(remember \ref{cl1}(a+b)), and thus $z_t\in I$. If $\alpha$ is limit of
uncountable cofinality the same conclusion is even more immediate.  

So $\{z_t:t\in T^+\}\subseteq I$, and $\bB/{\rm id}_\bB(\{z_t:t\in T^+\})$
is a Boolean Algebra generated by its set of atoms which is included in  
\[\{y\in {\mathcal A}_2:\neg(\exists t)(y\leq z_t)\}\]
(by \ref{cl1}(a)). Hence we conclude that $\bB$ is superatomic.
\end{proof}

\begin{claim}
\label{cl3}
\begin{enumerate}
\item $\Atom(\bB)=\{x_{(t,\alpha)}:(t,\alpha)\in T\times\mu\}$, so $\bB$ has
$|T|+\mu$ atoms, and $|\bB|=\chi$. 
\item $|\Aut(\bB)|\leq 2^\mu$; moreover for every $f\in\Aut(\bB)$
\[|\{(t,\alpha)\in T\times\mu:f(x_{(t,\alpha)})\ne x_{(t,\alpha)}\}|<\mu.\]
\end{enumerate}
\end{claim}

\begin{proof}[Proof of the Claim]
(1)\quad Easy.\\
(2)\quad Clearly the second statement implies the first. So let $f\in\Aut(
\bB)$ and suppose that $f$ moves at least $\mu$ atoms. Then there is $g\in
{\mathcal F}$ such that $f(x_{(t,\alpha)})=x_{g(t,\alpha)}$ for all $(t,
\alpha)\in\Dom(g)$. But, by \ref{cl1}(c), there is $y\in {\mathcal A}_2$
such that $y\subseteq\Dom(g)$ and $g[y]$ is almost disjoint to every member
of ${\mathcal A}_2$.  An easy contradiction. 
\end{proof}
\end{proof}

\begin{remark}
\label{2.2} 
\begin{enumerate}
\item As $\End(\bB)\ge|\bB|$ this gives an example for \cite[Problem
76]{M2}, too.  Still the first example works in more cardinals and is 
different.
\item With a little more work we can guarantee that the number of one-to-one endomorphisms of $\bB$ is $\leq 2^\mu$.
\item Alternatively, for the proof of \ref{cl1} we can use $\mu^{\aleph_0}=
2^\mu$ almost disjoint subsets of $\Dom(g_\alpha)$, say $\langle y_{\alpha,
i}: i<2^\mu\rangle$; for each $i$ choose $y'_{\alpha,i}\in [y_{\alpha,i}]^{
\aleph_0}$ such that it satisfies clause (b) of \ref{cl1} (exists by Ramsey
theorem), so for some $i$ we have: $y_{\alpha,i}$, $y'_{\alpha,i}=:g_\alpha[
y_{\alpha,i}]$ are almost disjoint to $y_\beta,y'_\beta$ for $\beta<\alpha$.
So $y_{\alpha,i},y'_{\alpha,i}$ are as required.
\end{enumerate}
\end{remark}

\section{On entangledness}

\begin{definition}
\label{3.3} 
\begin{enumerate}
\item A sequence $\bar{\mathcal I}=\langle {\mathcal I}_\alpha:\alpha<
\alpha^*\rangle$ of linear orders is $\kappa$-entangled if:
\begin{enumerate}
\item[(a)] each ${\mathcal I}_\alpha$ is a linear order of cardinality $\ge
\kappa$, and 
\item[(b)] if $n<\omega$, $\alpha_1<\ldots<\alpha_n<\alpha^*$, and
$t^\ell_\zeta\in {\mathcal I}_{\alpha_\ell}$ for $\ell\in \{1,\ldots,n\}$,
$\zeta<\kappa$ are such that $\zeta\neq\xi\ \Rightarrow\ t^\ell_\zeta\neq
t^\ell_\xi$, then for any $w\subseteq\{1,\ldots,n\}$ we may find $\zeta<\xi
<\kappa$ such that:
\[\ell\in w\ \Rightarrow\ {\mathcal I}_{\alpha_\ell}\models t^\ell_\zeta<
t^\ell_\xi\qquad\mbox{and}\qquad\ell\in\{1,\ldots,n\}\setminus w\
\Rightarrow\ {\mathcal I}_{\alpha_\ell}\models t^\ell_\xi<t^\ell_\zeta.\]
\end{enumerate}
If $\kappa$ is omitted we mean: $\kappa=\min\{|{\mathcal I}_\alpha|:\alpha<
\alpha^*\}$.  
\item $\Ens(\kappa,\lambda)$ is the statement asserting that there is an
entangled sequence $\bar{\mathcal I}=\langle {\mathcal I}_\alpha:\alpha<
\lambda\rangle$ of linear orders each of cardinality $\kappa$.
\end{enumerate}
\end{definition}

\begin{definition}
\label{3.2}
\begin{enumerate}
\item For an ideal $J$ on $\kappa$ we let
\[\begin{array}{ll}
\bU_J(\chi)=:\min\{|{\mathcal A}|:&{\mathcal A}\subseteq [\chi]^\kappa
 \mbox{ and}\\
&(\forall f\in {}^\kappa\chi)(\exists A\in {\mathcal A})(\{i<\kappa:f(i)\in
 A\}\in J^+)\}.
  \end{array}\]
\item $\Ded^+(\mu)=:\min\{\theta$: there is no linear order with $\theta$
elements and density $\leq\mu\;\}$. 
\end{enumerate}
\end{definition}

\begin{theorem}
\label{3.1} 
Assume that $\mu<\kappa<\chi<\Ded^+(\mu)$, $2^\mu<\lambda$, $\kappa$ is
regular and $\lambda\leq\bU_{J^{\bd}_\kappa}(\chi)$ (see Definition
\ref{3.2}). Then $\Ens(\kappa,\lambda)$ by a sequence $\langle {\mathcal
I}_\alpha:\alpha<\lambda\rangle$ of linear orders of cardinality $\kappa$
and density $\mu$ (see Definition \ref{3.3}).
\end{theorem}

\begin{proof}
Let ${\mathcal J}$ be a dense linear order of cardinality $\chi$ with a
dense subset ${\mathcal J}^*$ of cardinality $\mu$. Without loss of
generality the set of elements of ${\mathcal J}$ is $\chi$ and of ${\mathcal 
J}^*$ is $\mu$. Let $u^i_\zeta$ (for $i<\kappa$, $\zeta<\chi$) be pairwise
distinct members of ${\mathcal J}$, and let $\bar{u}=\langle u^i_\zeta:i<
\kappa,\zeta<\chi\rangle$. For $f\in {}^\kappa\chi$ let ${\mathcal I}_f=
\{u^i_{f(i)}:i<\kappa\}$. 

\begin{mclaim}
\label{3.5} 
If $n<\omega$, $f_0,\ldots,f_{n-1}\in {}^\kappa\chi$ and $\bar{\mathcal I}=
\langle {\mathcal I}_{f_\ell}:\ell<n\rangle$ is entangled, then we can find
${\mathcal A} \subseteq [\chi]^\kappa$ such that $|{\mathcal A}|\le 2^\mu$
and: 
\begin{enumerate}
\item[$(\oplus)$] if $f\in {}^\kappa\chi$ and $(\forall A\in {\mathcal A})(
\forall^* i<\kappa)(f(i)\notin A)$, then $\bar{\mathcal I}\conc\langle
{\mathcal I}_f\rangle$ is entangled 

(\/$\forall^*$ means ``for every large enough''). 
\end{enumerate}
\end{mclaim}

\begin{proof}[Proof of the Claim]
Assume $f_0,\ldots,f_{n-1}\in {}^\kappa\chi$ and $\bar{\mathcal I}=\langle
{\mathcal I}_{f_\ell}:\ell<n\rangle$ is entangled. 

Let 
\[{\mathcal F}=\{f\in {}^\kappa\chi:\bar{\mathcal I}\conc\langle
{\mathcal I}_f\rangle\mbox{ is not entangled}\;\}.\]
For each $f_n=f\in{\mathcal F}$ we fix $w^f\subseteq\{0,\ldots,n\}$ and
$t^{\ell,f}_j\in {\mathcal I}_{f_\ell}$ (for $\ell\leq n$, $j<\kappa$) with
no repetitions witnessing that $\bar{\mathcal I}\conc\langle{\mathcal I}_f
\rangle$ is not entangled. Next we fix a model $N_f\prec ({\mathcal H}(
\beth^+_8(\chi)),\in,<^*)$ such that $\mu+1\subseteq N_f$, $\|N_f\|=\mu$,
$\{\bar{\mathcal I},{\mathcal I}_f,{\mathcal J},f\}\in N_f$ and $\bar{t}^f= 
\langle t^{\ell,f}_j:\ell\leq n,j<\kappa\rangle\in N_f$. Note that for
$i<\kappa$ we have: 
\begin{enumerate}
\item[(i)] $t^{\ell,f}_i\notin N_f$ whenever $i\notin N_f$, 
\item[(ii)] $x\in N_f\ \&\ |x|<\kappa\ \&\ \sup(N_f\cap\kappa)\le i<\kappa
\quad\Rightarrow\quad i\notin x$.
\end{enumerate}

Now we define a relation $E$ on ${\mathcal F}$ letting for $f,g\in {\mathcal
F}$: 
\[\begin{array}{lcl}
f\;E\;g&\mbox{ if and only if }&(\alpha)\quad w^f=w^g,\\
& &(\beta)\quad N_f\cap\chi=N_g\cap\chi,\\
& &(\gamma)\quad (\forall\ell\leq n)(\forall j\in N_f\cap\kappa)(t^{\ell,
f}_j=t^{\ell,g}_j).
  \end{array}\]
Note that $E$ is an equivalence relation on ${\mathcal F}$, and there are at
most $2^\mu$ $E$--equivalence classes. Therefore, in order to show
\ref{3.5}, it is enough that for each $E$--equivalence class $g/E$ we define 
a set $Y_{g/E}\in [\chi]^\kappa$ such that:  
\begin{enumerate}
\item[$(\boxtimes)$]\qquad if $f\in g/E$ then $\neg (\forall^*i)(f(i)\notin
Y_{g/E})$. 
\end{enumerate}
Then, letting ${\mathcal A}=\{Y_{g/E}:g\in {\mathcal F}\}$ we will get a
family as required in \ref{3.5}.  
\medskip

So let $g\in{\mathcal F}$, $w^*=w^g$ and let $i^*=\sup(N_g\cap\kappa)$.

For $i<\kappa$, and a sequence $\bar{t}=\langle t^\ell:\ell<n\rangle\in
\prod\limits_{\ell<n}{\mathcal I}_{f_\ell}$ we let 
\[Y^i_{\bar{t}}=\big\{f(j): f\in g/E\ \&\ (\forall\ell<n)(t^{\ell,f}_i=
t^\ell)\ \&\ j<\kappa\ \&\ u^j_{f(j)}=t^{n,f}_i\}.\]
We claim that 
\begin{enumerate}
\item[(iii)] if $i>i^*$ (but $i<\kappa$) and $\bar{t}\in\prod\limits_{\ell<
n}{\mathcal I}_{f_\ell}$, then $|Y^i_{\bar{t}}|\leq 1$.
\end{enumerate}
Why? Assume toward contradiction that $f_1(j_1),f_2(j_2)$ are two distinct
members of $Y^i_{\bar{t}}$, $f_1,f_2\in g/E$, $t^{\ell,f_m}_i=t^\ell$ (for
$\ell<n$ and $m=1,2$) and $t^{n,f_1}_i=u^{j_1}_{f_1(j_1)}\neq u^{j_2}_{f_2(
j_2)}=t^{n,f_2}_i$. Pick disjoint intervals $(a^1,b^1), (a^2,b^2)$ of
${\mathcal J}^*$ such that $t^{n,f_m}_i\in (a^m,b^m)$ and $[t^\ell\neq
t^{n,f_m}_i\ \Rightarrow\ t^\ell\notin (a^m,b^m)]$ (for $m=1,2$ and
$\ell<n$). Without loss of generality, if $n\in w^*$ then $b^1<a^2$, else
$b^2<a^1$. We can also pick $a_\ell,b_\ell\in {\mathcal J}^*$ (for $\ell<n$)
such that $a_\ell<t^\ell<b_\ell$ and:   
\begin{itemize}
\item if $t^\ell\neq t^{\ell'}$ then $(a_\ell,b_\ell)\cap (a_{\ell'},
b_{\ell'})=\emptyset$,
\item if $t^\ell\neq t^{n,f_m}_i$ then $(a_\ell,b_\ell)\cap (a^m,b^m)=
\emptyset$.
\end{itemize}
Now, we are going to show that
\begin{enumerate}
\item[(iii)$^*$] if $w\subseteq n$, $m\in\{1,2\}$, and $i_0\in N_{f_m}\cap
\kappa$, and $a^+_\ell\in {\mathcal J}^* \cap [a_\ell,t^\ell)$ and
$b^+_\ell\in {\mathcal J}^*\cap (t^\ell,b_\ell]$ (for $\ell<n)$, 

\noindent then we can find $j\in N_{f_m}\cap\kappa\setminus i_0$ such that 
\begin{enumerate}
\item[$(*)_j$] $t^{n,f_m}_j\in (a^m,b^m)$, and 
\[\ell\in w\ \Rightarrow\ a^+_\ell<t^{\ell,f_m}_j<t^\ell\qquad\mbox{ and
}\qquad\ell\in n\setminus w\ \Rightarrow\ t^\ell<t^{\ell,f_m}_j<b^+_\ell.\] 
\end{enumerate}
\end{enumerate}
So assume that (iii)$^*$ fails, so there is no $j\in N_m\cap\kappa\setminus
i_0$ such that $(*)_j$ holds. First note that then also there is no $j'<i$
(but $j'>i_0$) satisfying $(*)_{j'}$. [Why? Suppose $(*)_{j'}$ holds true
and choose $a^*_\ell,b^*_\ell\in {\mathcal J}^*$ such that 
\[\begin{array}{rcl}
\ell\in w& \Rightarrow& a^+_\ell=a^*_\ell<t^{\ell,f_m}_{j'}<b^*_\ell<t^\ell
\qquad\mbox{ and}\\
\ell\in n\setminus w&\Rightarrow& t^\ell<a^*_\ell<t^{\ell,f_m}_{j'}<b^*_\ell
=b^+_\ell.
  \end{array}\] 
The set 
\[Z=\{j\in\kappa\setminus i_0:(\forall\ell<n)(a^*_\ell<t^{\ell,f_m}_j<
b^*_\ell)\ \&\ a_m<t^{\ell,f_m}_j<b_m\}\]
is non-empty (as witnessed by $j'$) and it belongs to the model $N_{f_m}$. 
Picking any $j'\in Z\cap N_{f_m}$ provides a witness for (iii)$^*$ (so we
get a contradiction).]

\noindent Next, the set 
\[Z_0=:\{j<\kappa:(\forall\ell<n)(t^{\ell,f_m}_j\in (a^+_\ell,b^+_\ell))\
\&\ t^{n,f_m}_j\in (a^m,b^m)\}\] 
belongs to $N_{f_m}$ and $i$ belongs to it. But $i>i^*$, so necessarily
$Z_0$ has cardinality $\kappa$ (remember (ii)). Let 
\[Z_1=:\{j\in Z_0\setminus i_0:(\exists j_1<j)(j_1\in Z_0\ \&\ (\forall\ell
<n)(t^{\ell,f_m}_{j_1}<t^{\ell,f_m}_j\equiv\ell\in w))\}.\]
By the assumption that (iii)$^*$ fails (and the discussion above) we have
$i\notin Z_1$. But again $Z_1 
\in N_{f_m}$, so $|Z_0\setminus Z_1|=\kappa$. Since the sequence
$\bar{\mathcal I}$ is entangled, we can find $j_1<j_2$ in $Z_0\setminus Z_1$
such that $(\forall\ell<n)(t^{\ell,f_m}_{j_1}<t^{\ell,f_m}_{j_2}\equiv
\ell\in w)$. But then $j_1$ witnesses $j_2\in Z_1$, a contradiction. 
\medskip

Now we are going to use (iii)$^*$ twice to justify (iii). First we apply
(iii)$^*$ for $w=:w^*$, $i_0=0$, $m=1$ with $a^+_\ell=a_\ell$, $b^+_\ell=
b_\ell$ getting $j_1\in N_{f_1}\cap\kappa$ such that  
\[\begin{array}{rcl}
\ell\in w^*& \Rightarrow& a_\ell<t^{\ell,f_1}_{j_1}<t^\ell,\qquad\mbox{
and}\\  
&& \\
\ell\in n\setminus w^*&\Rightarrow& t^\ell<t^{\ell,f_1}_{j_1}<b_\ell,
  \end{array}\]
and $t^{n,f_1}_{j_1}\in (a^1,b^1)$. Next we choose $a^+_\ell,b^+_\ell\in
{\mathcal J}^*$ (for $\ell<n$) such that
\[\begin{array}{rcll}
\ell\in w^*&\Rightarrow&t^{\ell,f_1}_{j_1}<a^+_\ell<t^\ell&\mbox{and }\
b^+_\ell=b_\ell,\\ 
\ell\in n\setminus w^*&\Rightarrow&t^\ell<b^+_\ell<t^{\ell,f_1}_{j_1}&
\mbox{and }\ a^+_\ell=a_\ell.
  \end{array}\]
Then we again apply (iii)$^*$, this time for $w=:w^*$, $m=2$, $i_0=j_1+1$
and $a^+_\ell,b^+_\ell$ chosen above, getting $j_2\in N_{f_2}\cap\kappa
\setminus j_1$ such that, in particular, $(\forall\ell<n)(t^{\ell,f_2}_{j_2}
\in(a^+_\ell,b^+_\ell))$ and $t^{\ell,f_2}_{j_2}\in (a^2,b^2)$. Then clearly 
\[(\forall\ell\le n)(t^{\ell,f_1}_{j_1}<t^{\ell,f_2}_{j_2}\equiv \ell\in
w^*),\]
and $j_1<j_2$ both are in $N_{f_1}\cap\kappa=N_{f_2}\cap\kappa$. Since
$f_1,f_2$ are $E$-equivalent we know that $t^{\ell,f_1}_{j_1}=t^{\ell,
f_2}_{j_1}$ (for $\ell\leq n$), so we may get a contradiction with the
choice of $\bar{t}^{f_2}$ and we finish the proof of (iii). 
\smallskip

Now we let 
\[Y_{g/E}=\bigcup\{Y^i_{\bar{t}}: i^*<i<\kappa\ \&\ \bar{t}\in
\prod_{\ell<n}{\mathcal I}_{f_\ell}\}.\]
It follows from (iii) that $|Y_{g/E}|\leq\kappa$. Clearly, for each $f\in
g/E$ the set $\{j<\kappa: f(j)\in Y_{g/E}\}$ is of size $\kappa$. Hence
$Y_{g/E}$ is as required in $(\boxtimes)$ and this finishes the proof of
\ref{3.5}.  
\end{proof}

\noindent{\bf Continuation of the proof of \ref{3.1}:}\qquad Now we can
construct the entangled sequence of linear orders as required in the
theorem. For this, by induction on $\alpha<\lambda$, we choose functions
$f_\alpha\in {}^\kappa\chi$ such that: 
\begin{enumerate}
\item[$(\otimes_\alpha)$]\qquad the sequence $\langle{\mathcal I}_{f_\beta}:
\beta<\alpha\rangle$ is entangled.
\end{enumerate}
Note that if $\alpha\leq\lambda$ is limit and $f_\beta$ have been chosen for
$\beta<\alpha$ so that $(\otimes_\beta)$ holds (for $\beta<\alpha$), then
also $(\otimes_\alpha)$ holds. Let $f_0\in {}^\kappa\chi$ be any function;
note that $(\otimes_1)$ holds true as $\kappa$ is $>\mu$ which is the
density of ${\mathcal J}$, so in ${\mathcal J}$ there is no monotonic
sequence of length $\mu^+$.

Suppose we have defined $f_\beta\in{}^\kappa\chi$ for $\beta<\alpha$ so that
$(\otimes_\alpha)$ holds true. Let $\langle\bar{\beta}^\zeta:\zeta<\alpha^*
\rangle$ list all the sequences $\langle\beta_\ell:\ell<n\rangle\subseteq
\alpha$ such that $n<\omega$ and $\bigwedge\limits_{\ell_1\neq\ell_2}
\beta_{\ell_1}\neq\beta_{\ell_2}$. Let $\bar{\beta}^\zeta=\langle\beta(
\zeta,\ell):\ell<n_\zeta\rangle$. Clearly without loss of generality 
\[\alpha^*=|\alpha|\ \vee\ (\alpha<\omega\ \&\ \alpha^*<\omega).\] 
For each $\zeta<\alpha^*$ we apply \ref{3.5} to $f_{\beta(\zeta,0)},\ldots,
f_{\beta(\zeta,n_\zeta-1)}$ to get a family ${\mathcal A}^\zeta \subseteq
[\chi]^\kappa$ as there (so in particular $|{\mathcal A}^\zeta|\le 2^\mu$). 
There is $f_\alpha\in {}^\kappa\chi$ such that 
\[(\forall\zeta<\alpha^*)(\forall A\in {\mathcal A}^\zeta)(\forall^* i< \kappa)
(f_\alpha(i)\notin A).\]
Why? Otherwise $\bigcup\limits_{\zeta<\alpha^*}{\mathcal A}^\zeta$
exemplifies 
\[\bU_{J^{\bd}_\kappa}(\chi)\leq |\bigcup\limits_{\zeta<\alpha^*} {\mathcal
A}^\zeta|\leq (|\alpha|\!+\!\aleph_0)\cdot \sup\{|{\mathcal A}^\zeta|:
\zeta\!<\!\alpha^*\}\leq (|\alpha|\!+\!\aleph_0)\times 2^\mu<\lambda.\]
Now, with $f_\alpha$ chosen as above, $(\otimes_{\alpha+1})$ holds true. 
\end{proof}

\begin{remark}
Theorem \ref{3.1} should be compared with:
\begin{enumerate}
\item[(a)] \cite[Ch.II, 4.10E]{Sh:g}, see AP2 there on history. There
we got only $\Ens_2$. 
\item[(b)] \cite[\S2]{Sh:462}, but there the density is higher. 
\end{enumerate}
\end{remark}

\begin{conclusion}
\label{3.6}   
\begin{enumerate}
\item Let $\kappa$ be an uncountable regular cardinal $\le 2^{\aleph_0}$,
$\kappa<\chi\leq 2^{\aleph_0}$, and $\bU_{J^{\bd}_\kappa}(\chi)>
2^{\aleph_0}$ (e.g., $\chi=2^{\aleph_0}$, $\cf(\chi)=\kappa<\chi$). Then
there is an entangled sequence of length $\bU_{J^{\bd}_\kappa}(\chi)$ of
linear orders of cardinality $\kappa$.  
\item Assume $\mu$ is a strong limit singular cardinal, $\mu<\kappa= \cf(
\kappa)<\chi\leq 2^\mu$ and $\bU_{J^{\bd}_\kappa}(\chi)>2^\mu$ (e.g., $\chi=
2^\mu$, $\cf(\chi)=\kappa<\chi$). Then there is an entangled sequence of
length $\bU_{J^{\bd}_\kappa}(\chi)$ of linear orders of cardinality $\kappa$.
\end{enumerate}
\end{conclusion}

\section{On attainment of spread}
In this section we are interested in the following question 

\begin{question}
\label{5.0} 
Let $\lambda$ be a singular cardinal.
\begin{enumerate}
\item Is there a Boolean algebra $\bB$ such that $s^+(\bB)=\lambda$, e.g.,
in the following sense: 
\begin{quotation}
there is no sequence $\langle a_\alpha:\alpha<\lambda\rangle\subseteq\bB
\setminus\{0\}$ such that each $a_\alpha$ is not in the ideal generated by
\[I_\alpha=\{a_\beta:\beta\neq\alpha\},\]
but for each $\mu<\lambda$ there is such a sequence? 
\end{quotation}
\item We can ask also/alternatively for $\hd^+(\bB)=\lambda$ (and/or
$\hL^+(\bB)=\lambda$) defined similarly using $\{a_\beta:\beta<\alpha\}$
(and/or $\{a_\beta:\beta>\alpha\}$, respectively). 
\end{enumerate}
\end{question}
For the discussion of the attainment properties of spread we refer the
reader to \cite[p.~175]{M2}; the attainment of $\hd$, $\hL$ is discussed,
e.g., in \cite[p.~198, p.~191]{M2}. Forcing constructions for different 
attainment properties for $\hd$ and $\hL$ are presented in \cite{RoSh:651}.

\begin{theorem}
\label{5.1}
\begin{enumerate}
\item Assume that $\mu$ is a strong limit singular cardinal, 
\[\aleph_0<\cf(\mu)<\mu<\cf(\lambda)<\lambda \le 2^\mu.\]
Then
\begin{enumerate}
\item[$(\boxtimes_\lambda)$] there is a Boolean Algebra $\bB$ satisfying:
\begin{enumerate}
\item[(i)]   $|\bB|=\lambda=s(\bB)$,
\item[(ii)]  $s(\bB)$ is not obtained (i.e., $s^+(\bB)=\lambda$),
\item[(iii)] moreover $\hd^+(\bB)=\hL^+(\bB)=\lambda$.
\end{enumerate}
\end{enumerate}
\item Assume that
\begin{enumerate}
\item[$(\otimes_2)$] 
\begin{enumerate}
\item[(a)] $\mu<\cf(\lambda)<\lambda$, 
\item[(b)] $\langle\lambda_i:i<\delta\rangle$ is a (strictly) increasing
sequence of regular cardinals with limit $\mu$, 
\item[(c)] $J$ is an ideal on $\delta$ extending $J^{\bd}_\delta$, $A\in
J^+$, $\delta\setminus A \in J^+$, 
\item[(d)] $\langle g_\alpha:\alpha<\cf(\lambda)\rangle$ is a $<_{J
\restriction A}$--increasing $<_{J\restriction A}$--cofinal sequence of
members of $\prod\limits_{i\in A}\lambda_i$, and $\langle h_\alpha:\alpha<
\lambda\rangle$ is a sequence of distinct members of $\prod\limits_{i\in
\delta\setminus A}\lambda_i$ such that  
\[j<\delta\quad \Rightarrow\quad |\{h_\alpha\restriction j,g_\beta
\restriction j:\alpha<\lambda,\ \beta<\cf(\lambda)\}|<\lambda_j.\]
\end{enumerate}
\end{enumerate}
Then $(\boxtimes_\lambda)$ holds.
\item Assume that
\begin{enumerate}
\item[$(\otimes_3)$]
\begin{enumerate}
\item[(a)] $\mu<\cf(\lambda)<\lambda$,
\item[(b)] $\langle\lambda_i:i<\delta\rangle$ is a strictly increasing
sequence of regular cardinals $<\mu$, 
\item[(c)] $J$ is an ideal on $\delta$ extending $J^{\bd}_\delta$, $A
\subseteq\delta$, $A\in J^+$ and $\delta\setminus A\in J^+$,
\item[(d)] $g_\alpha\in\prod\limits_{i<\delta}\lambda_i$ for $\alpha<
\lambda$ are pairwise distinct, 
\item[(e)] among $\{g_\alpha\restriction A:\alpha<\lambda\}$ we can find an
$<_{J\restriction A}$--increasing cofinal sequence of length $\cf(\lambda)$, 
\item[(f)] $|\{g_\alpha\restriction i:\alpha<\lambda\}|=\lambda_i$.
\end{enumerate}
\end{enumerate}
Then $(\boxtimes_\lambda)$ holds.
\end{enumerate}
\end{theorem}

\begin{proof}  1)\quad We shall prove that the assumptions of part (2)
hold. 

As $\cf(\mu)>\aleph_0$, we know (by \cite[Ch.VIII, \S1]{Sh:g}) that there is
a sequence $\langle\lambda_i:i<\cf(\mu)\rangle$ such that 
\[\mu>\lambda_i=\cf(\lambda_i)>|\prod\limits_{j<i}\lambda_j|\quad\mbox{ and
}\quad\tcf(\prod\limits_{i<\cf(\mu)}\lambda_i/J^{\bd}_{\cf(\mu)})=
\cf(\lambda).\]
Let $\langle g_\alpha:\alpha<\cf(\lambda)\rangle$ be an increasing cofinal
sequence in $(\prod\limits_{i<\cf(\mu)}\lambda_i,<_{J^{\bd}_{\cf(\mu)}})$. 
Let $h_\alpha\in\prod\limits_i\{\lambda_{2i+1}:i<\kappa\}$ (for $\alpha\in
[\cf(\lambda),\lambda)$) be just such that $h_\alpha\notin\{h_\beta:\beta<
\alpha\}$, so $A=:\{2i:i<\cf(\mu)\}$, $\langle g_\alpha\restriction A:
\alpha<\cf(\lambda)\rangle$, $\langle h_\alpha \restriction (\kappa\setminus
A): \alpha<\lambda\rangle$ are as required $(\otimes_2)$. 
\medskip

\noindent 2)\quad Let $\langle\chi_i:i<\cf(\lambda)\rangle$ be an increasing
continuous sequence of cardinals such that 
\begin{itemize}
\item $\lambda=\sum\limits_{i<\cf(\lambda)}\chi_i$, 
\item $\chi_0=0$, $\cf(\lambda)<\chi_1$ and each $\chi_{i+1}$ is regular. 
\end{itemize}
For $\alpha<\lambda$ let $j(\alpha)<\cf(\lambda)$ be such that $\alpha\in 
[\chi_{j(\alpha)},\chi_{j(\alpha)+1})$ and let $f_\alpha\in\prod\limits_{i<
\delta}\lambda_i$ be such that: 
\[f_\alpha\restriction A=g_{j(\alpha)}\qquad\mbox{ and }\qquad f_\alpha
\restriction (\delta\setminus A)=h_\alpha.\]

Now for $n\geq 1$ we define a Boolean Algebra $\bB_n$ (each $\bB_n$ will be
an example):

it is generated by $\{x_\alpha:\alpha<\lambda\}$ freely except:
\begin{enumerate}
\item[$(\circledast)$] \quad {\em if\/} $i\in A$, $\nu_k\in\prod\limits_{i'
<i}\lambda_{i'}$, $\nu_k\conc\langle\gamma_{k,\ell}\rangle\vartriangleleft
f_{\alpha_{k,\ell}}$ (for $k<m$, $\ell\leq 2n+1$), and $w\subseteq m$, and 
\[\begin{array}{rcl}
\ell< 2n\ \&\ k\in m\setminus w&\Rightarrow&\gamma_{k,\ell}<\gamma_{k,\ell+
1},\\
k\in w&\Rightarrow &\alpha_{k,n}=\alpha_{k,2n+1},
  \end{array}\]
and there are no repetitions in the sequence $\langle\nu_k:k<m\rangle$, and
$\bt_k\in\{0,1\}$, 

\noindent {\em then}\qquad $\bigcap\limits_{k<m} x^{\bt_k}_{\alpha_{k,n}}
\leq\bigcup\limits_{\moj{\scriptstyle\ell\ne n,}{\scriptstyle\ell\leq 2n+1}}
\bigcap\limits_{k<m} x^{\bt_k}_{\alpha_{k,\ell}}$,  
\end{enumerate}
where $x^\bt$ is $x$ if $\bt=1$, and $-x$ if $\bt=0$.

\begin{claim}
\label{Fact}
$s^+(\bB_n)\leq\lambda$, $\hd^+(\bB_n)\leq\lambda$, $\hL^+(\bB_n)\leq
\lambda$. 
\end{claim}

\begin{proof}[Proof of the Claim]
Assume toward contradiction that the sequence $\langle a_\beta:\beta<\lambda
\rangle\subseteq\bB_n\setminus\{0\}$ exemplifies the failure. Without loss
of generality, $a_\beta=\bigcap\limits_{\ell<m_\beta}x^{\bt(\beta,
\ell)}_{\alpha(\beta,\ell)}$, where $\ell<m<m_\beta\ \Rightarrow\ \alpha(
\beta,\ell)\ne\alpha(\beta,m)$. For each $i<\cf(\lambda)$ we choose $S_i
\subseteq [\chi_i,\chi_{i+1})$, and $\varepsilon_i(*)<\delta$, $m^i<\omega$,
$\bt[i,\ell]\in\{0,1\}$, $j[i,\ell]<\cf(\lambda)$ (for $\ell<m^i$) such that
(note that we can permute $\langle\alpha(\beta,\ell):\ell<m_\beta\rangle$):  
\begin{enumerate}
\item[(i)]   $S_i$ is unbounded in $\chi_{i+1}$,
\item[(ii)]  for all $\beta\in S_i$ we have
\[m_\beta=m^i\ \&\ (\forall\ell<m^i)(\bt(\beta,\ell)=\bt[i,\ell]\ \&\
j(\alpha(\beta,\ell))=j[i,\ell]),\] 
\item[(iii)] $\langle\langle\alpha(\beta,\ell):\ell<m^i\rangle:\beta\in S_i
\rangle$ is a $\Delta$--system with heart $\langle\alpha[i,\ell]:\ell<k^i
\rangle$, so 
\[\begin{array}{rcl}
\beta\in S_i\ \&\ \ell<k^i& \Rightarrow&\alpha(\beta,\ell)=\alpha[i,\ell],
 \qquad\mbox{ and}\\
\alpha(\beta_1,\ell_1)=\alpha(\beta_2,\ell_2)&\Rightarrow&(\beta_1=\beta_2 \
 \&\ \ell_1=\ell_2)\vee (\ell_1=\ell_2< k^i),
  \end{array}\]
\item[(iv)]  for $\beta\in S_i$, there are no repetitions in the sequence
$\langle f_{\alpha(\beta,\ell)}\restriction\varepsilon_i(*):\ell<m^i\rangle$
and it does not depend on $\beta$,
\item[(v)]   for every $\beta^*\in S_i$ and $\varepsilon<\delta$ the set 
\[\{\beta\in S_i:(\forall\ell<m^i)(f_{\alpha(\beta,\ell)}\restriction
\varepsilon=f_{\alpha(\beta^*,\ell)}\restriction \varepsilon)\}\]
is unbounded in $S_i$.
\end{enumerate}
Note that necessarily 
\begin{enumerate}
\item[(vi)]  $j[i,\ell]\ge i$ for $\ell\in [k^i,m^i)$.
\end{enumerate}
Next pick a set $S\in [\cf(\lambda)]^{\cf(\lambda)}$ such that:
\begin{enumerate}
\item[$(\alpha)$] for all $i\in S$ we have $m^i=m^*$, $k^i=k^*$, $\bt[i,
\ell]= \bt[\ell]$, $\varepsilon_i(*)=\varepsilon(*)$, 
\item[$(\beta)$]  $\langle\langle\alpha[i,\ell]:\ell<k^*\rangle:i\in S
\rangle$ is a $\Delta$--system with heart $\langle\alpha(\ell):\ell<\ell^*
\rangle$, so 
\[i\in S\ \&\ \ell<\ell^*\quad\Rightarrow\quad\alpha[i,\ell]=\alpha(\ell),\]
\item[$(\gamma)$] also $\langle\langle j[i,\ell]:\ell<m^*\rangle:i\in S
\rangle$ is a $\Delta$--system with heart $\langle j(\ell):\ell\in w^*
\rangle$, where $w^*\subseteq m^*$. 
\end{enumerate}
Note that then $\ell^*\subseteq w^*\subseteq k^*$ (the first inclusion is a
consequence of $(\beta)$, the second one follows from (vi)). 

Also by further shrinking of the sets $S_i$ (for $i<\cf(\lambda)$) and $S$
we may require that 
\begin{enumerate}
\item[(A)] if $i_1<i_2$ are from $S$, then $j[i_1,\ell]<i_2$ (for
$\ell<m^*$),  
\item[(B)] if $i_1\ne i_2$ are from $S$ and $\beta_1\in S_{i_1}$ and
$\beta_2\in S_{i_2}$, then 
\[\{\alpha(\beta_1,\ell):\ell<m^*\}\cap\{\alpha(\beta_2,\ell):\ell<m^*\}
\subseteq\{\alpha(\ell):\ell<\ell^*\},\]
\item[(C)] if $i_1\in S$, $\gamma_1\in S_{i_1}$, then 
\[(\forall\xi<\delta)(\exists^{\cf(\lambda)}i\in S)(\exists^{\chi_{i+1}}
\gamma\in S_i)(\forall\ell<m^*)(f_{\alpha(\gamma,\ell)}\restriction\xi =
f_{\alpha(\gamma_1,\ell)}\restriction \xi).\] 
\end{enumerate}
Choose $\gamma_i\in S_i$ for $i\in S$. Look at $\bar{f}^i=\langle f_{\alpha(
\gamma_i,\ell)}:\ell<m^*\rangle$.

We can (as in \cite[Ch.II, 4.10A]{Sh:g}) find $\varepsilon<\delta$ and
$\bar{f}=\langle f_0,\ldots,f_{m^*-1}\rangle$ such that $\varepsilon\in A$,
$\varepsilon>\varepsilon(*)$ and: 
\begin{enumerate}
\item[$(*)$]  for every $\zeta<\lambda_\varepsilon$ there is $i\in S$ such
that: 
\[(\forall\ell<m^*)\big(f_{\alpha(\gamma_i,\ell)}\restriction\varepsilon=
f_\ell\restriction\varepsilon\big)\ \mbox{ and }\ (\forall\ell\in m^*
\setminus w^*)\big(f_{\alpha(\gamma_i,\ell)}(\varepsilon)>\zeta\big).\]
\end{enumerate}
So we can choose inductively $\zeta_k,i_k$ (for $k\leq 2n$) such that
$i_k\in S$, $\zeta_k<\lambda_\varepsilon$, and
\[(\forall\ell\!<\!m^*)\big(f_{\alpha(\gamma_{i_k},\ell)}\!\restriction\!
\varepsilon= f_\ell\!\restriction\!\varepsilon\big)\ \mbox{ and }\ (\forall
\ell\!\in\! m^*\!\setminus\! w^*)\big(\zeta_k<f_{\alpha(\gamma_{i_k},\ell)}
(\varepsilon)<\zeta_{k+1}\big).\] 
Note that, as $\varepsilon \in A$, we have
\[(\forall\ell\in w^*)\big(f_{\alpha(\gamma_{i_k},\ell)}(\varepsilon)=
g_{j(\alpha(\gamma_{i_k},\ell))}(\varepsilon)=g_{j(\ell)}(\varepsilon)\big)\] 
for each $k\leq 2n$. It follows from clause (v) above that we may pick
$\gamma\in S_{i_n}\setminus\{\gamma_{i_n}\}$ such that $(\forall\ell<m^*)(
f_{\alpha(\gamma,\ell)}\restriction\varepsilon=f_{\alpha(\gamma_{i_n},\ell)} 
\restriction\varepsilon)$. By our choices, $\alpha(\gamma_{i_n},\ell)=\alpha
(\gamma,\ell)$ for $\ell<k^*$ (so in particular for $\ell\in w^*$). Now, by
the definition of $\bB_n$, we clearly have $a_{\beta_n}\leq
\bigcup\limits_{\moj{\scriptstyle\ell\neq n,}{\scriptstyle\ell\leq 2n+1}}
a_{\beta_\ell}$, where $\beta_\ell=\gamma_{i_\ell}$ for $\ell\leq 2n$ and
$\beta_{2n+1}=\gamma$, finishing the proof for $s$.  

Now for $\hd,\hL$ use clause (C) above.
\end{proof}

\begin{claim}
$s^+(\bB_n)>\chi_{i+1}$, more specifically $\{x_\alpha:\alpha\in [\chi_i,
\chi_{i+1})\}$ are independent as ideal generators.
\end{claim}

\begin{proof}[Proof of the Claim]
Let $\alpha^*\in [\chi_i,\chi_{i+1})$. We define a function $h_{\alpha^*}=h:
\{x_\alpha:\alpha<\lambda\}\longrightarrow\{0,1\}$ by:
\[h(x_\alpha)=\left\{
\begin{array}{ll}
1&\mbox{ if }\alpha=\alpha^*,\\
0&\mbox{ if }\llg(f_\alpha\cap f_{\alpha^*})\in\delta\setminus A,\\
1&\mbox{ if }\llg(f_\alpha\cap f_{\alpha^*})\in A,\ f_\alpha(\llg(f_\alpha
\cap f_{\alpha^*}))>f_{\alpha^*}(\llg(f_\alpha\cap f_{\alpha^*})),\\ 
0&\mbox{ if }\llg(f_\alpha\cap f_{\alpha^*})\in A,\ f_\alpha(\llg(f_\alpha
\cap f_{\alpha^*}))<f_{\alpha^*}(\llg(f_\alpha\cap f_{\alpha^*})).
\end{array}\right.\]
We claim that the function $h$ respects the equations in the definition of
$\bB_n$. To show this suppose that $i\in A$, $\bt_k\in\{0,1\}$,
$\nu_k\in\prod\limits_{i'<i}\lambda_{i'}$, $\nu_k\vartriangleleft
f_{\alpha_{k,\ell}}$ (for $k<m$, $\ell\leq 2n+1$) and $w\subseteq m$ are as
in the assumptions of $(\circledast)$. Now we consider three cases.
\smallskip

\noindent{\sc Case 1:}\qquad $f_{\alpha^*}\restriction i\notin\{\nu_k:
k<m\}$.\\
Then, by the way $h$ is defined, $h(x_{\alpha_{k,n}})=
h(x_{\alpha_{k,\ell}})$ for each $\ell\leq 2n+1$ and $k<m$. Hence easily
\[\bigcap_{k<m}h(x_{\alpha_{k,n}})^{\bt_k}=\bigcup_{\moj{\scriptstyle\ell
\neq n,}{\scriptstyle\ell\leq 2n+1}}\bigcap_{k<m} h(x_{\alpha_{k,\ell}}
)^{\bt_k},\]
and we are done.
\smallskip

\noindent{\sc Case 2:}\qquad $f_{\alpha^*}\restriction i=\nu_{k^*}$, $k^*\in 
m\setminus w$.\\
Thus $f_{\alpha_{k^*,0}}(i)<\ldots<f_{\alpha_{k^*,n}}(i)<\ldots<
f_{\alpha_{k^*,2n}}(i)$ and $h(x_{\alpha_{k^*,0}})=h(x_{\alpha_{k^*,n}})$ or
$h(x_{\alpha_{k^*,2n}})=h(x_{\alpha_{k^*,n}})$. Let $\ell^*$ be 0 in the
first case and $2n$ in the second. Note that also for $k<m$, $k\neq k^*$
we have $h(x_{\alpha_{k,\ell}})=h(x_{\alpha_{k,n}})$ for all $\ell\leq
2n+1$. Hence
\[\bigcap_{k<m}h(x_{\alpha_{k,\ell^*}})^{\bt_k}=\bigcap_{\moj{\scriptstyle
\ k\neq k^*,}{\scriptstyle k<m}}h(x_{\alpha_{k,n}})^{\bt_k}\cap
h(x_{\alpha_{k^*,\ell^*}})^{\bt_{k^*}}=\bigcap_{k<m}h(x_{\alpha_{k,n}})^{
\bt_k},\]
and we are done.  
\smallskip

\noindent{\sc Case 3:}\qquad $f_{\alpha^*}\restriction i=\nu_{k^*}$, $k^*\in 
w$.\\
Thus $\alpha_{k^*,n}=\alpha_{k^*,2n+1}$ (so $h(x_{\alpha_{k^*,n}})=h(
x_{\alpha_{k^*,2n+1}})$) and also for $k<m$, $k\neq k^*$ we have
$h(x_{\alpha_{k,n}})=h(x_{\alpha_{k,2n+1}})$. Hence
\[\bigcap_{k<m}h(x_{\alpha_{k,n}})^{\bt_k}=\bigcap_{k<m}
h(x_{\alpha_{k,2n+1}})^{\bt_k},\]
and we are done.
\smallskip

Consequently the function $h$ can be extended to a homomorphism $\hat{h}$
from $\bB_n$ to $\{0,1\}$. Clearly $h(x_{\alpha^*})=1$ and  $h(x_\alpha)=0$
for all $\alpha\in [\chi_i,\chi_{i+1})\setminus\{\alpha^*\}$. (Remember
$f_\alpha\restriction A=g_i$ for $\alpha\in [\chi_i,\chi_{i+1})$, and hence  
\[\mbox{if }\ \alpha\ne\beta\in [\chi_i,\chi_{i+1})\ \mbox{ then }\
\llg(f_\alpha\cap f_\beta)\in\delta\setminus A.\ \mbox{ )}\] 
Thus we are done.
\end{proof}
\medskip

\noindent 3)\quad We can get the assumptions of part (2). 
\end{proof}

\begin{remark}
\label{5.1A}
\begin{enumerate}
\item We cannot really prove in ZFC that there is a Boolean Algebra $\bB$
such that $s^+(\bB)$ is singular ($\equiv s(\bB)$ singular not obtained) as
$s^+(\bB)$ cannot be strong limit singular. 
\item Note that the demand $(\exists \mu)[\mu<\cf(\lambda)<\lambda<2^\mu]$
is necessary by \cite{Sh:233}. The construction is like the one in 
\cite[\S 7]{RoSh:599}. Earlier see \cite[4.14]{Sh:345}.
\item Of course, the proof of \ref{5.1}(2) shows that we have the respective
result for finite variants $s_m$ of spread, as well as for $\hd_m,\hL_m$ (if
$m\geq 3$, i.e., $m=2n+1$). We refer the reader to \cite[\S 1]{RoSh:534} for
the definitions of these cardinal invariants (see also \cite{RoSh:599} for
discussion and some independence results on $s_m$; more relevant results can
be found in \cite{Sh:620}).  
\end{enumerate}
\end{remark}

So we can give examples to \ref{5.0} if we can have for
$(\boxtimes_\lambda)$ of \ref{5.1}.

\begin{proposition}
\label{5.2} 
\begin{enumerate}
\item If $\kappa$ is strong limit singular cardinal, $2^\kappa\geq
\aleph_{\kappa^+}$, then we have examples of $\lambda$, $\kappa<\cf(\lambda)
<\lambda\leq 2^\kappa$ with $(\boxtimes_\lambda)$ (of \ref{5.1}), e.g.,
$\lambda=\aleph_{\kappa^+}$ !
\item If $\delta(*)=(2^\kappa)^{+(\kappa^{+4})}$, $\aleph_{\delta(*)}\le 
2^{\kappa^+}$, then also there is $\lambda\in (2^\kappa,2^{\kappa^+})$,
$\cf(\lambda)\in [2^\kappa,(2^\kappa)^{{+\kappa}^{+4}})$ as needed in
\ref{5.1}(3), and hence $(\boxtimes_\lambda)$.
\item If $\kappa$ is inaccessible (possibly weakly) $\delta(*)=
(2^{<\kappa})^{{+ \kappa}^{+4}}$ and $\aleph_{\delta(*)}<2^\kappa$ then we
can find $\lambda\in [2^{<\kappa},2^\kappa)$, $\cf(\lambda)\in [2^{<\kappa},
(2^{<\kappa})^{{+\kappa}^{+4}})$, as in \ref{5.1}(2), and hence
$(\boxtimes_\lambda)$ holds.\\
Similarly if $\kappa$ is a singular cardinal, or a successor cardinal by
part (2).  
\item E.g., if $\aleph_{\aleph_{\omega+1}}\le 2^{\aleph_0}$, then for
$\lambda =\aleph_{\aleph_{\omega+1}}$ we have an example for this cardinal.   
Generally, if $\mu>\cf(\mu)=\aleph_0$, $\cf(\lambda)=\mu^+$ and $\lambda \le
2^{\aleph_0}$, then there is an example in $\lambda$. 
\end{enumerate}
\end{proposition}

\begin{proof}
1)\quad Should be clear. (Note that $\pp^+(\kappa)\geq\aleph_{\kappa^++1}$
by \cite[5.9, p.408]{Sh:400}).
\medskip

\noindent 2)\quad First, for some club $C$ of $\kappa^{+4}$ (for $\alpha\in
C\ \Rightarrow\ \alpha$ limit) we have 
\begin{enumerate}
\item[$(*)_1$]\qquad $\delta\in C\ \&\ \cf(\delta)\le\kappa\quad\Rightarrow
\quad\pp((2^\kappa)^{+ \delta})<(2^\kappa)^{+\min(C\setminus (\delta+1))}$.
\end{enumerate}
(By \cite[\S4]{Sh:400}).  Hence (again by \cite{Sh:g}) 
\begin{enumerate}
\item[$(*)_2$]\qquad $\delta\in C\ \&\ \cf(\delta)\leq\kappa\quad\Rightarrow 
\quad((2^\kappa)^{+\delta})^\kappa<(2^\kappa)^{+\min(C\setminus(\delta+1))}$. 
\end{enumerate}
We can for any $\delta\in\acc(C)$ with $\cf(\delta)=\kappa^+$ do the
following: we can find a strictly increasing sequence $\langle\lambda_i:i<
\kappa^+\rangle$ of regular cardinals with limit $(2^\kappa)^{+\delta}$,
$2^\kappa<\lambda_i$, $\tcf(\prod\limits_{i<\kappa^+}\lambda_i)/
J^{\bd}_{\kappa^+}= (2^\kappa)^{+\delta +1}$ (if we assume
$\pp((2^\kappa)^{+ \delta})>(2^\kappa)^{+\delta +1}$ we can find more
examples). 

Note: 
\[j<\kappa^+\quad\Rightarrow\quad\prod_{i<j}\lambda_i<(2^\kappa)^{+\delta}\] 
(by $(*)_2$), as the ideal is $J^{\bd}_{\kappa^+}$; without loss of
generality 
\begin{enumerate}
\item[$(*)_3$]\qquad $\prod\limits_{i<j}\lambda_i<\lambda_j$. 
\end{enumerate}
So let $\langle g'_\alpha:\alpha<(2^\kappa)^{+\delta+1}\rangle$ be
$<_{J^{\bd}_{\kappa^+}}$--increasing and cofinal in $\prod\limits_{i<
\kappa^+}\lambda_i$ and let $A=\{2i:i<\kappa^+\}$. 

Now assume $2^{\kappa^+}\geq\lambda>\cf(\lambda)=(2^\kappa)^{+\delta +1}$;
such $\lambda$ exists by the assumption. We can find $h_\alpha\in
{}^{\kappa^+}2$ (for $\alpha\in [(2^\kappa)^{+\delta+1},\lambda))$ with no
repetitions. 

Note $|\{g'_\alpha\restriction i,h'_\alpha\restriction i:\alpha<\lambda\}|
\le |\prod\limits_{j<i}\lambda_j|$, which has cardinality $<\lambda_i$. So
we can apply Theorem \ref{5.1}(3). 
\medskip 

\medskip 3), 4)\quad  Same. 
\end{proof}

\begin{discussion}
\label{5.3}
\begin{enumerate}
\item If Cardinal Arithmetic is too close to GCH ($2^\kappa<\aleph_{
\kappa^+}$ for every $\kappa$), no example exists as by \cite{Sh:233}, ${\rm
ZFC}\models 2^{\cf(s^+(\bB))}>|\bB|$. [Why? If $\bB$ is a counterexample,
let $\lambda=s^+(\bB)=s(\bB)$ (bring a counterexample); clearly $\lambda$ is
a limit cardinal, so $2^{\cf(\lambda)}>|\bB|\geq\lambda\geq\aleph_{\cf(
\lambda)}$, a contradiction.] 

If Cardinal Arithmetic is far enough from GCH (even just for regulars), then
there is an example. 

I consider it a semi-ZFC answer --- see \cite{Sh:666} and \cite{Sh:430}. 

\item There are some variants of problem \ref{5.0} related to various
versions of the (equivalent) definitions of $s,\hd,\hL$. For $s$ all
versions are equivalent \cite[p. 175]{M2}. Concerning $\hd,\hL$ see
the discussion of the attainment relations for the equivalent
definitions of $\hd$ in \cite[p. 196, 197]{M2} and of $\hL$ in
\cite[p.191]{M2}. On the remaining cases see also in \cite[\S 4]{RoSh:651}.  

\end{enumerate}
\end{discussion}

\begin{problem}
Does $\aleph_{\omega_1}<2^{\aleph_0}$ imply that an example for $\lambda=
\aleph_{\omega_1}$ exists? 
\end{problem}

\bibliographystyle{bplain}
\bibliography{listb,lista,liste,listx}

\def\germ{\frak} \def\scr{\cal}
  \ifx\documentclass\undefinedcs\def\rm{\fam0\tenrm}\fi
  \def\defaultdefine#1#2{\expandafter\ifx\csname#1\endcsname\relax
  \expandafter\def\csname#1\endcsname{#2}\fi} \defaultdefine{Bbb}{\bf}
  \defaultdefine{frak}{\bf} \defaultdefine{mathbb}{\bf}
  \defaultdefine{mathcal}{\cal}
  \defaultdefine{beth}{BETH}\defaultdefine{cal}{\bf} \def\bbfI{{\Bbb I}}
  \def\mbox{\hbox} \def\text{\hbox} \def\om{\omega} \def\Cal#1{{\bf #1}}
  \def\pcf{pcf} \defaultdefine{cf}{cf} \defaultdefine{reals}{{\Bbb R}}
  \defaultdefine{real}{{\Bbb R}} \def\restriction{{|}} \def\club{CLUB}
  \def\w{\omega} \def\exist{\exists} \def\se{{\germ se}} \def\bb{{\bf b}}
  \def\equivalence{\equiv} \let\lt< \let\gt> \def\cite#1{[#1]}
\begin{thebibliography}{10}
\makeatletter \renewcommand{\@biblabel}[1]{[#1]} \makeatother

\bibitem{M2}
J.~Donald Monk.
\newblock {\em {Cardinal Invariants of Boolean Algebras}}, volume 142 of {\em
  {Progress in Mathematics}}.
\newblock Birkh\"auser Verlag, Basel--Boston--Berlin, 1996.


\bibitem{RoSh:651}
Andrzej Roslanowski and Saharon Shelah.
\newblock {Forcing for hL and hd}.
\newblock {\em {Colloquium Mathematicum}}, submitted.
{\bf [RoSh:651]};   {\tt math.LO/9808104}\footnote{References of the
form {\tt math.XX/$\cdots$} refer to the {\tt xxx.lanl.gov} archive}


\bibitem{RoSh:534}
Andrzej Roslanowski and Saharon Shelah.
\newblock {Cardinal invariants of ultrapoducts of Boolean algebras}.
\newblock {\em {Fundamenta Mathematicae}}, 155:101--151, 1998.
{\bf [RoSh:534]};   {\tt math.LO/9703218}

\bibitem{RoSh:599}
Andrzej Roslanowski and Saharon Shelah.
\newblock {More on cardinal functions on Boolean algebras}.
\newblock {\em {Annals of Pure and Applied Logic}}, 103:1--37, 2000.
{\bf [RoSh:599]};   {\tt math.LO/9808056}

\bibitem{Sh:E12}
Saharon Shelah.
\newblock {Analytical Guide and Corrections to \cite{Sh:g}.}
{\bf [Sh:E12]};   {\tt math.LO/9906022}

\bibitem{Sh:666}
Saharon Shelah.
\newblock {On what I do not understand (and have something to say)}.
\newblock {\em {Fundamenta Mathematicae}}, to appear.
{\bf [Sh:666]};   {\tt math.LO/9906113}

\bibitem{Sh:233}
Saharon Shelah.
\newblock {Remarks on the numbers of ideals of Boolean algebra and open sets of
  a topology}.
\newblock In {\em {Around classification theory of models}}, volume 1182 of
  {\em {Lecture Notes in Mathematics}}, pages 151--187. {Springer, Berlin},
  1986.
{\bf [Sh:233]}

\bibitem{Sh:345}
Saharon Shelah.
\newblock {Products of regular cardinals and cardinal invariants of products of
  Boolean algebras}.
\newblock {\em {Israel Journal of Mathematics}}, 70:129--187, 1990.
{\bf [Sh:345]}

\bibitem{Sh:g}
Saharon Shelah.
\newblock {\em {Cardinal Arithmetic}}, volume~29 of {\em {Oxford Logic
  Guides}}.
\newblock {Oxford University Press}, 1994.
{\bf [Sh:g]}

\bibitem{Sh:400}
Saharon Shelah.
\newblock {Cardinal Arithmetic}.
\newblock In {\em {Cardinal Arithmetic}}, volume~29 of {\em {Oxford Logic
  Guides}}, chapter~IX. {Oxford University Press}, 1994.
\newblock Note: See also [Sh400a] below.
{\bf [Sh:400]}

\bibitem{Sh:430}
Saharon Shelah.
\newblock {Further cardinal arithmetic}.
\newblock {\em {Israel Journal of Mathematics}}, 95:61--114, 1996.
{\bf [Sh:430]};   {\tt math.LO/9610226}

\bibitem{Sh:462}
Saharon Shelah.
\newblock {On $\sigma $-entangled linear orders}.
\newblock {\em {Fundamenta Mathematicae}}, 153:199--275, 1997.
{\bf [Sh:462]};   {\tt math.LO/9609216}

\bibitem{Sh:620}
Saharon Shelah.
\newblock {Special Subsets of ${}^{{\rm cf}(\mu)}\mu$, Boolean Algebras and
  Maharam measure Algebras}.
\newblock {\em {Topology and its Applications}}, 99:135--235, 1999.
\newblock 8th Prague Topological Symposium on General Topology and its
  Relations to Modern Analysis and Algebra, Part II (1996).
{\bf [Sh:620]};   {\tt math.LO/9804156}

\end{thebibliography}
\end{document}